%% file: orgdt_tps.tex
\documentclass{IEEEtran}

\ifCLASSINFOpdf
\else
\fi

\usepackage{amsmath}
\usepackage{psfrag}
\usepackage{multicol}
\usepackage{epstopdf}
\usepackage{graphicx}
\PassOptionsToPackage{bookmarks=false}{hyperref}
\usepackage{subfigure,gensymb,url}
\usepackage{color}
\usepackage{anysize, booktabs, setspace}
\marginsize{0.8in}{0.8in}{0.8in}{0.8in}
\usepackage{amsfonts, amssymb, amsthm}
\usepackage{cleveref, nomencl, blindtext, mathtools, etoolbox}
\crefname{section}{§}{§§}   
\Crefname{section}{§}{§§}
\usepackage{pgfplots,tikz}

\usepackage[ruled,vlined,linesnumbered]{algorithm2e}

\renewcommand\nomgroup[1]{%
  \item[\bfseries
  \ifstrequal{#1}{A}{Parameters and notation}{%
  \ifstrequal{#1}{B}{Binary variables}{%
  \ifstrequal{#1}{C}{Continuous variables}{%
  \ifstrequal{#1}{D}{Parameters}{%
  \ifstrequal{#1}{E}{Random variables}{%
  \ifstrequal{#1}{F}{Scenario parameters}{%
   \ifstrequal{#1}{X}{Other Symbols}{}}}}}}}]%
  }

\immediate\write18{makeindex \jobname.nlo -s nomencl.ist -o \jobname.nls}
\makenomenclature
\newcommand{\ORGDT}{{\mathcal{P}_0^{AC}}}

\newcommand{\SubProblem}{{\cal{Q}^{AC}}(s)}
\newcommand{\SubProblemDC}{{\cal{Q}^{DC}}(s)}

\newcommand{\Nodes}{{\cal{N}}}
\newcommand{\Edges}{{\cal{E}}}
\newcommand{\Disasters}{{\cal{S}}}
\newcommand{\Damages}{{\cal{D}}}

\newcommand{\resist}{{R}}
\newcommand{\react}{{X}}

\newcommand{\lineVariable}{x}
\newcommand{\switchVariable}{\tau}
\newcommand{\hardenVariable}{t}

\newcommand{\facilityVariable}{u}

\newcommand{\loss}{{l^s}}

\newcommand{\capacity}{T}

\newcommand{\Load}{{\cal{L}}}

\allowdisplaybreaks

\newtheorem{thm}{Theorem}[section]

\newtheorem{rmk}[thm]{Remark}
\newtheorem{obs}[thm]{Observation}

\begin{document}
\title{Resilient Transmission Grid Design: AC Relaxation vs. DC approximation}

\author{
   Harsha Nagarajan$^{\dag}$, Russell Bent$^{\dag}$, Pascal Van Hentenryck$^\ddag$, Scott Backhaus$^{\dag}$, Emre Yamangil$^\dag$ \\
$^\dag$ Center for Nonlinear Studies, Los Alamos National Laboratory, NM, United States. \\
$^\ddag$ Department of Industrial Operations Engineering, University of Michigan, Ann Arbor, United States.\\
Contact: harsha@lanl.gov
}

\maketitle

\begin{abstract}

As illustrated in recent years (Superstorm Sandy, the Northeast Ice Storm of 1998, etc.), extreme weather events pose an enormous threat to the electric power transmission systems and the associated socio-economic systems that depend on reliable delivery of electric power. Besides inevitable malfunction of power grid components, deliberate malicious attacks can cause high risks to the service. These threats motivate the need for approaches and methods that improve the resilience of power systems. In this paper, we develop a model and tractable methods for optimizing the upgrade of transmission systems through a combination of hardening existing components, adding redundant lines, switches, generators, and FACTS and phase-shifting devices. While many of these controllable components are included in traditional design (expansion planning) problems, we uniquely assess their benefits from a resiliency point of view. 
More importantly, perhaps, we evaluate the suitability of using state-of-the-art AC power flow relaxations versus the common DC approximation in resilience improvement studies. The resiliency model and algorithms are tested on a modified version of the RTS-96 (single area) system.

\end{abstract}

 \input{Intro}

\input{Formulation}

 \input{Algo}
\input{NR}

 \input{Conclusions}

\bibliographystyle{IEEEtran}
\bibliography{references.bib}

\end{document}

%% file: Intro.tex

\section{Introduction}
\label{Sec:intro}


\IEEEPARstart{T}{he} modern electrical system is designed for transportation of large amounts of power from sources of supply to distant points of demand. Within these systems, the underlying high-voltage transmission networks play a vital role in achieving this mission. However, when transmission networks are exposed to extreme event conditions, the ability to deliver power is degraded because of physical damage to overhead 
transmission lines and towers. One example of such events are ice storms. During an ice storm, transmission towers can fail due to leg buckling
and lines can fail due to the combined stress of ice accumulation and wind \cite{Brostrom2007,Eidinger,wang2010influence}. 

When such events occur on large scales, outages and impacts can be extreme. For example, in the winter of 1998, an ice storm in northeastern North America toppled over 1000 transmission towers and 30,000 wooden utility poles. Over 5 million people were without power and the economic impacts were estimated at \$2.6 billion \cite{Mills2012}. Thus, given the potential social and economic impacts of these events, it is important to consider how to upgrade the design of transmission systems to improve their performance under these conditions.

In our preliminary work \cite{7540988}, we formulated the Optimal Resilient Grid Design problem for Transmission systems (ORGDT) 
as a two-stage mixed-integer stochastic optimization problem and developed an algorithm 
to solve this problem. The first stage selects from a set of potential upgrades to the network. The
second stage evaluates the network performance benefit of the upgrades 
with a set of damage scenarios sampled from a stochastic distribution of events of concern. 
One of the important questions that was unanswered in this early work centered on the use of recently developed convex relaxations of AC power flow physics \cite{hijazi2014convex} in resiliency planning. Here, we address this concern and convincingly show that using DC approximation of power flow physics leads to severely suboptimal and infeasible solutions under certain operating conditions of the grid. These results establish the necessity of using models of AC physics in resilient design problems and analysis.  

Like \cite{7540988}, we adopt the methods discussed in \cite{wang2010influence} in order to sample 
realistic damage scenarios for transmission systems. The ORGDT upgrade options include: a) Build new lines, b) Build switches, FACTS devices and Phase-shifting transformers (PST) to provide operational flexibility, c) Harden existing lines to lower the 
probability of damage, and d) Build new generation facilities.
Minimal network (resiliency) performance is measured by satisfying a minimum fraction of critical and non-critical loads. We use the exact algorithms developed in \cite{7540988} to compare solutions found using the DC approximation and the convex quadratic relaxations of \cite{hijazi2014convex}.


\noindent \textbf{Literature Review} 
Our recent paper \cite{7540988} is the most closely related work. This paper develops the ORGDT model and algorithm used here and is based on the resilient distribution system design work developed in \cite{yamangil2015resilient}. 

Another important area of related work is interdiction modeling and optimization. Here, the
goal is to operate or design a system to make it as resilient as possible to an adversary who can damage up to $k$ elements \cite{Chen:13,Chen2014a,Salmeron2009,Delgadillo2010}. These models are a generalization of our model when $k$ is chosen to bound a worst-case event.
However, given their min-max structure, these models are computationally challenging and are solvable only for small $k$. Instead, we exploit the
probabilistic nature of the adversary and we are able to address larger problems. In addition, existing interdiction models also do not generally include AC physics. 
Closely related to interdiction and our model is reference \cite{nezamoddini2017risk}.  They develop a model for hardening and upgrading transmission systems that is resilient to all possible physical damages bounded by a certain budget (i.e. weighted $k$). Instead of developing an attacker-defender model as in interdiction modeling, 
they include each possible damage as a scenario (like our model). Given the combinatorial number of possible scenarios, they develop heuristics for reducing the number of scenarios that are included. Unlike our approach, they focus only on the DC approximation and leave the development of decomposition algorithms for future work.
A third area of related work is stochastic transmission and generation expansion planning, where a recent survey describes some of the state-of-the-art \cite{Krishnan2015}. 

Overall, most papers in network design and expansion planning use the linearized DC model and few studies consider FACTS devices and PSTs, although they may have significant benefits. Some notable exceptions include the use of PSTs in network expansion \cite{miasaki2012transmission}, which uses
a genetic algorithm over the DC model. See also the recent work in \cite{Frolov2014} for the optimal placement of these devices to avoid congestion. 

The key contributions of this paper include:

\begin{itemize}
\item A detailed case study and extensions of the ORGDT model (from \cite{7540988}) that compares DC approximation based solutions and solutions based on state-of-the-art AC power flow relaxations. 
\item An approach for recovering AC feasible solutions from solutions found using the DC approximation and AC relaxations.
\item A detailed expansion model with appropriate generator capacity modeling, FACTS devices and PST devices, and cutting-plane-based algorithmic enhancements that was omitted from \cite{7540988}.
\end{itemize}

The rest of this paper is organized as follows. In section \ref{Sec:formulation} we describe the ORGDT model. Section \ref{Sec:algo} describes our decomposition algorithms and Section \ref{Sec:case} our case studies. We conclude with Section \ref{Sec:conc}.

%% file: Formulation.tex
\section{ORGDT Optimization Model}
\label{Sec:formulation}

{\fontsize{8.3}{6}\selectfont
\input{Nomenclature}
\printnomenclature
}

We formulate the ORGDT as a two-stage mixed-integer nonlinear program, $\ORGDT$. The first-stage variables (without superscript $s$) specify the new infrastructure enhancements and the second-stage variables (with superscript $s$) describe the operation of invested infrastructure to serve loads for each damage scenario $s \in \Disasters$. 

{\fontsize{8.7}{8}\selectfont
\begin{subequations}
\begin{flalign}
\ORGDT := \min & \sum_{ij\in \Edges} \left( c^{x}_{ij} x_{ij} + c^{\tau}_{ij} \tau_{ij} + c^{t}_{ij} t_{ij} + c^{\delta}_{ij}\delta_{ij} + c^{\gamma}_{ij}\gamma_{ij}\right)  \nonumber\\
& + \sum_{i \in \Nodes}\left(c^{u}_{i} u_{i} + c^{zp}_{i} zp_{i} \right)  \label{eqn:obj}  \\
\text{s.t.}\quad & x^s_{ij} \leq x_{ij} \quad \quad \forall ij\in \Damages_s, s\in \Disasters  \label{eqn:assign:1a} \\
& x^s_{ij} = x_{ij} \quad \quad \forall ij\notin \Damages_s, s\in \Disasters  \label{eqn:assign:1b} \\
& \tau^s_{ij} \leq \tau_{ij}, \  t^s_{ij} \leq t_{ij}   \quad \quad \forall ij\in \Edges, s\in \Disasters  \label{eqn:assign:2} \\
& \delta^s_{ij} \leq \delta_{ij}, \  \gamma^s_{ij} \leq \gamma_{ij}   \quad \quad \forall ij\in \Edges, s\in \Disasters  \label{eqn:assign:3} \\
& zp^s_{i} \leq zp_{i}, \ zq^s_{i} \leq zq_{i} \quad \forall i\in \Nodes, s\in \Disasters  \label{eqn:assign:4a} \\
& 0 \leq zp_{i} \leq zp^u_{i} u_{i}, \ |zq_{i}| \leq zq^u_{i} u_{i} \quad \forall i\in \Nodes   \label{eqn:assign:6a} \\
& zp_{i} \geq 2|zq_{i}|\quad \forall i\in \Nodes   \label{eqn:assign:6b} \\
& (\mathbf{x}^s,\boldsymbol{\tau}^s,\mathbf{t}^s,\boldsymbol{\delta}^s, \boldsymbol{\gamma}^s, \mathbf{zp}^s, \mathbf{zq}^s,\mathbf{u}) \in \SubProblem \quad \forall s\in \Disasters \label{eqn:feasible_network} \\
& x_{ij},\tau_{ij}, t_{ij},\delta_{ij},\gamma_{ij},u_{i} \in \{0,1\} \quad   \forall ij\in \Edges, i\in \Nodes \label{eqn:discrete_2}
\end{flalign}
\label{eq:master_AC}
\end{subequations}}
\noindent
In $\ORGDT$, Eq. \eqref{eqn:obj} minimizes the total upgrade cost, which includes new lines, switches, hardening, FACTS, phase shifters, generators and real power capacity, respectively. Constraints \eqref{eqn:assign:1a}-\eqref{eqn:assign:4a} link the first-stage (construction) decisions with second-stage variables in $\SubProblem$.
Constraints \eqref{eqn:assign:6a} and \eqref{eqn:assign:6b} model generation 
capacity upgrades. Without loss of generality, reactive power capacity is limited to half of the real power capacity.
Constraint \eqref{eqn:assign:6b} ensures that the generators are not built purely for reactive power support. Constraint \eqref{eqn:feasible_network} states that the mixed-integer vector $(\mathbf{x}^s,\boldsymbol{\tau}^s,\mathbf{t}^s,\boldsymbol{\delta}^s, \boldsymbol{\gamma}^s, \mathbf{zp}^s, \mathbf{zq}^s,\mathbf{u}) \in \SubProblem$ is an AC feasible transmission network for scenario $s$ subject to the constraints of $\SubProblem$. 

The constraints of $\SubProblem$ describe the AC-power flow equations and budget constraints on resiliency options. 
Ohm's law is modeled in constraints \eqref{eqn:acopf_real_ij}-\eqref{eqn:acopf_reac_ji}. For a given topology of the network. $\tilde{x}^s$, flow on line $ij$ is forced to zero when $\tilde{x}^s_{ij}=0$. 
Kirchhoff's current law (flow conservation) is given in constraints \eqref{eqn:balance1} and \eqref{eqn:balance2}.
Constraints \eqref{eqn:acopf_loss1}-\eqref{eqn:acopf_loss3} model the connection between current magnitude ($l^s_{ij}$) and the power loss equations as discussed in \cite{baran1989optimal}.
Operational constraints \eqref{eqn:thermalij}-\eqref{eqn:v_bnd} represent the thermal limits, phase angle difference limits 
and voltage bounds, respectively. Constraint \eqref{eqn:damage} models the damaged lines of the scenario $s\in S$, i.e., a line is inoperable when damaged and unhardened. Constraint \eqref{eqn:lineminusswitch} defines the topology for scenario $s$. A switch is available for operation only if the line is active. Constraints \eqref{eqn:mg_a} bound the real and reactive power capacities that can be built (if $u_i=1$) at bus $i\in \Nodes$. Constraints in \eqref{eqn:gp_a} and \eqref{eqn:gq_a} represent total real and reactive power capacities including the existing generators of the network. Finally, constraint \eqref{eqn:shed} captures the continuous power shedding subject to the resiliency constraints in \eqref{eqn:p_cr}-\eqref{eqn:q_ncr}. 
User-defined resiliency parameters ($lp_{cr},lq_{cr},lp_{ncr},lq_{ncr}$) ensure that a minimum fraction of critical and non-critical loads is served during every damage scenario $s \in \Disasters$. 
\vspace{-.48cm}

{\fontsize{8.7}{8}\selectfont
\begin{flalign}
\SubProblem  &  = \{\mathbf{x}^s,\boldsymbol{\tau}^s,\mathbf{t}^s,\boldsymbol{\delta}^s, \boldsymbol{\gamma}^s, \mathbf{zp}^s, \mathbf{zq}^s,\mathbf{u} :  \nonumber \\
& \textbf{On/Off AC power flow equations:} \nonumber \\
& p^s_{ij} = \tilde{x}^s_{ij} (\tilde{G}_{ij} {v^s}^{2}_{i} - \tilde{G}_{ij}v^s_i v^s_j \cos(\theta^s_{ij}) \nonumber \\  
& -\tilde{B}_{ij}v^s_i v^s_j \sin(\theta^s_{ij}))  \quad  \forall ij\in \Edges,  \label{eqn:acopf_real_ij} \\
&q^s_{ij} = \tilde{x}^s_{ij} (-\tilde{B}_{ij} {v^s}^{2}_{i} + \tilde{B}_{ij}v^s_i v^s_j \cos(\theta^s_{ij})\nonumber \\ 
& - \tilde{G}_{ij}v^s_i v^s_j \sin(\theta^s_{ij}) )  \quad  \forall ij\in \Edges,  \label{eqn:acopf_reac_ij} \\
&p^s_{ji} = \tilde{x}^s_{ij} (\tilde{G}_{ij} {v^s}^{2}_{j} - \tilde{G}_{ij}v^s_i v^s_j \cos(\theta^s_{ij}) \nonumber \\  
& +\tilde{B}_{ij}v^s_i v^s_j \sin(\theta^s_{ij}) ) \quad  \forall ij\in \Edges,  \label{eqn:acopf_real_ji} \\
&q^s_{ji} = \tilde{x}^s_{ij} (-\tilde{B}_{ij} {v^s}^{2}_{j} + \tilde{B}_{ij}v^s_i v^s_j \cos(\theta^s_{ij})\nonumber \\ 
& + \tilde{G}_{ij}v^s_i v^s_j \sin(\theta^s_{ij}) ) \quad  \forall ij\in \Edges,  \label{eqn:acopf_reac_ji} \\
& \theta^s_{ij} =  \theta^s_{i} -  \theta^s_{j} + \phi^s_{ij} \quad  \forall ij\in \Edges, \label{eqn:ph_ang}\\
& \tilde{G}_{ij} = (\bar{G}_{ij} - G_{ij})\delta^s_{ij} + G_{ij} \quad \quad \forall ij \in \Edges \label{eqn:conduct_disj}, \\
& \tilde{B}_{ij} = (\bar{B}_{ij} - B_{ij})\delta^s_{ij} + B_{ij} \quad \quad \forall ij \in \Edges \label{eqn:suscept_disj}, \\
& p^s_{ij} + p^s_{ji} = \resist_{ij} \loss_{ij} \quad  \forall ij\in \Edges,  \label{eqn:acopf_loss1} \\
& q^s_{ij} + q^s_{ji} = \react_{ij} \loss_{ij} \quad  \forall ij\in \Edges,  \label{eqn:acopf_loss2} \\
& {p^s}^{2}_{ij} + {q^s}^{2}_{ij} = \loss_{ij} {v^s}^{2}_{i} \quad  \forall ij\in \Edges,  \label{eqn:acopf_loss3} \\
& gp^s_{i} - lp^s_{i} = \sum_{ij\in \Edges} p^s_{ij} + \sum_{ji\in \Edges} p^s_{ij} \quad \quad \forall i\in \Nodes,  \label{eqn:balance1}\\
& gq^s_{i} - lq^s_{i} = \sum_{ij\in \Edges} q^s_{ij} + \sum_{ji\in \Edges} q^s_{ij} \quad \quad \forall i\in \Nodes, \label{eqn:balance2}\\
& \textbf{Operational limits and topology constraints:} \nonumber \\
& {p^s_{ij}}^2+{q^s_{ij}}^2 \leq \tilde{x}^s_{ij} \capacity^2_{ij} \quad \quad \forall ij\in \Edges,  \label{eqn:thermalij} \\  
& {p^s_{ji}}^2+{q^s_{ji}}^2 \leq \tilde{x}^s_{ij} \capacity^2_{ij} \quad \quad \forall ij\in \Edges,  \label{eqn:thermalji} \\
& |\theta^s_{ij}| \leq \tilde{x}^s_{ij}\theta^u + (1-\tilde{x}^s_{ij})\theta^M \quad \quad  \forall ij \in \Edges, \label{eqn:theta_bnd}\\
& v^l_i \leq v^s_i \leq v^u \quad  \forall i\in \Nodes,  \label{eqn:v_bnd} \\
& -\phi^u \gamma^s_{ij} \leq \phi^s_{ij} \leq \phi^u \gamma^s_{ij} \quad \forall ij\in \Edges,  \label{eqn:phi_bnd} \\
& x^s_{ij} = t^s_{ij} \quad \quad \forall ij\in \Damages_s  \label{eqn:damage}\\
& \tilde{x}^s_{ij} = x^s_{ij} - \tau^s_{ij}\geq 0 \quad \quad \forall ij\in \Edges  \label{eqn:lineminusswitch}\\
& \delta^s_{ij} \leq \tilde{x}^s_{ij}, \ \gamma^s_{ij} \leq \tilde{x}^s_{ij} \quad \quad \forall ij \in \Edges, \\
& \textbf{Generation constraints}\quad  \forall i \in \Nodes: \nonumber \\
& 0 \leq zp^s_{i} \leq zp^u_{i}u_{i}, \ |zq^s_{i}| \leq zq^u_{i}u_{i}, \ zp^s_{i} \geq 2|zq^s_{i}| \label{eqn:mg_a}\\
& 0 \leq gp^s_{i} \leq  gp^u_{i} + zp^s_{i},  \label{eqn:gp_a}\\
& (gq^l_{i} - |zq^s_{i}|) \leq gq^s_{i} \leq (gq^u_{i} + |zq^s_{i}|),  \label{eqn:gq_a}\\
& \textbf{Resilience constraints:} \nonumber \\
& lp^s_{i} =  yp^s_{i} dp_{i}, \ lq^s_{i} =  yq^s_{i} dq_{i} \quad \quad \forall i\in \Nodes,  \label{eqn:shed}\\
& \sum_{i \in \Load} lp^s_{i} \geq lp_{cr} \sum_{i \in \Load} dp_{i}, \label{eqn:p_cr}\\
& \sum_{i \in \Nodes \setminus \Load} lp^s_{i} \geq lp_{ncr} \sum_{i \in \Nodes \setminus \Load} dp_{i} \label{eqn:p_ncr}\\
& \sum_{i \in \Load} lq^s_{i} \geq lq_{cr} \sum_{i \in \Load} dq_{i} \label{eqn:q_cr}\\
& \sum_{i \in \Nodes \setminus \Load} lq^s_{i} \geq lq_{ncr} \sum_{i \in \Nodes \setminus \Load} dq_{i} \label{eqn:q_ncr}\\
& \mathbf{x}^s,\boldsymbol{\tau}^s,\mathbf{t}^s,\boldsymbol{\delta}^s, \boldsymbol{\gamma}^s \in \{0,1\}; \quad 0 \leq yp^s, yq^s \leq 1 \nonumber \}
\end{flalign}
}
\vspace{-1cm}

\subsection{Additional Investment Options}
\label{subsec:facts}
In this paper, we also include FACTS and PST devices and  as investment options. These devices are often useful for addressing congestion in overloaded transmission systems, a situation that can occur during major disruptions.
In addition, these devices are cost-effective  \cite{Baldick2014} and may reduce resiliency costs significantly by replacing the need for new transmission lines or hardening the existing damaged lines. 

\textbf{Series Compensators (FACTS)} We model series compensation with reactance reduction. As described in the nomenclature, $\delta^s_{ij}$ indicates if a compensation device is used on line $ij$ during scenario $s \in \Disasters$. To model this disjunction we introduce $(\tilde{G}_{ij}, \tilde{B}_{ij})$, a tuple representing the following conditional expression: 
\[
    (\tilde{G}_{ij},\tilde{B}_{ij})= 
\begin{cases}
    (\bar{G}_{ij},\bar{B}_{ij}),& \text{if used } (\delta^s_{ij}=1)\\       
    (G_{ij}, B_{ij}), & \text{otherwise } (\delta^s_{ij}=0)
\end{cases}
\]
This expression is modeled as constraints  \eqref{eqn:conduct_disj} and \eqref{eqn:suscept_disj}.
Further, this disjunction is captured in the modified power flow equations \eqref{eqn:acopf_real_ij}-\eqref{eqn:acopf_reac_ji}. 
The non-convexity of this disjunction is discussed later in this section. 

\textbf{Phase-Shifting Transformers (PSTs)} PSTs are devices that adjust voltage phase angles. 
As described in the nomenclature, $\gamma^s_{ij}$ indicates if a phase shifter
is used on line $ij$ during scenario $s \in \Disasters$. This disjunction is described in the following conditional expression:
\[
    \theta^s_{ij}= 
\begin{cases}
    \begin{aligned}
    &\theta^s_{i}-\theta^s_{j}+\phi^s_{ij}, \\
    &-\phi^u \leq \phi^s_{ij} \leq \phi^u.
    \end{aligned} & \text{if installed } (\gamma^s_{ij}=1), \\
    \theta^s_{i}-\theta^s_{j}, &\text{otherwise } (\gamma^s_{ij}=0).
\end{cases}
\]
This expression is modeled through constraints  \eqref{eqn:ph_ang} and \eqref{eqn:phi_bnd}. 


\subsection{Convex Relaxations for On/Off AC Power Flow}
\label{subsec:convex}
Numerous convex relaxations of AC power flow equations have been proposed in the literature to obtain tight lower bounds on the original nonconvex formulation. These relaxations have various trade-offs in complexity, quality, and performance. However, very few relaxation methods can be applied in the context of transmission switching problems as the addition of binary variables increases the complexity of the problems significantly. Hence, given the mixed-integer nature of ORGDT, which leads to an even more complex formulation with various resiliency options, we extend the recently developed convex quadratic relaxations of the power flow equations to ORGDT \cite{hijazi2014convex}. For completeness, we discuss the necessary convex relaxations for multilinear and trigonometric expressions which appear in nonconvex constraints \eqref{eqn:acopf_real_ij}-\eqref{eqn:acopf_reac_ji} and \eqref{eqn:acopf_loss3}.

\label{subsec:qc}
\noindent
\textbf{Multilinear expressions} 
Given any two variables $x_i$, $x_j \in \left[\underline{x}_i \overline{x}_i\right]\times \left[\underline{x}_j\overline{x}_j\right]$, the  McCormick relaxation is used to linearize the bilinear product $x_ix_j$ by introducing a new variable $\widehat{x_{ij}} \in {\langle x_i, x_j \rangle}^{MC}$. The feasible region of this variable is defined by equations (\ref{eq:SMC}). This relaxation is exact if one variable is binary \cite{nagarajan2016tightening}. 
\vspace{-0.5cm}

{\fontsize{9}{8}\selectfont
\begin{subequations} \label{eq:SMC}
\small
\begin{align}
& \label{McCormick}\widehat{x_{ij}} \geq \underline{x}_ix_j+\underline{x}_jx_i -\underline{x}_i \hspace{2pt} \underline{x}_j \\ 
&\widehat{x_{ij}} \geq \overline{x}_ix_j+\overline{x}_jx_i - \overline{x}_i \hspace{2pt} \overline{x}_j \\ 
&\widehat{x_{ij}} \leq \underline{x}_ix_j+\overline{x}_jx_i-\underline{x}_i \hspace{2pt} \overline{x}_j \\ 
&\widehat{x_{ij}} \leq \overline{x}_ix_j+\underline{x}_jx_i-\overline{x}_i \hspace{2pt} \underline{x}_j
\end{align}
\end{subequations}%
}
\noindent
Multilinear expressions of the form $x_ix_j\ldots x_p$ are relaxed using a sequential bilinear approach as follows: $\langle \langle x_i,x_j  \rangle^{MC},\ldots, x_p\rangle^{MC}$. 

\noindent \textbf{Quadratic terms}
Given a variable $x_i \in \left[\underline{x}_i \overline{x}_i\right]$, a second-order conic relaxation is applied to $x_i^2$ by introducing a new variable $\widehat{x_{i}} \in {\langle x_i\rangle}^{MC-q}$, such that,
\vspace{-.25cm}
\begin{subequations} 
\label{eq:MC-q}
\small
\allowdisplaybreaks
\begin{align}
& \label{McC-q}\widehat{x_{i}} \geq x_i^2\\ 
&\widehat{x_{i}} \leq (\overline{x}_i+\underline{x}_i)x_i - \overline{x}_i\underline{x}_i
\end{align}
\end{subequations}%

By applying the above approaches and disjunctive convex relaxations for the trigonometric functions with on/off variables \cite{hijazi2014convex}, we outer-approximate  
the non-convex power flow constraints. For brevity, we focus on constraint \eqref{eqn:acopf_real_ij} but similar relaxations hold  for constraints \eqref{eqn:acopf_reac_ij}-\eqref{eqn:acopf_reac_ji}.

\vspace{-0.2cm}
{\fontsize{9}{8}\selectfont
\begin{align*}
p^s_{ij} = &(\overline{G}_{ij} - G_{ij})\left( \widehat{\delta xv}^s_{ij} - \widehat{\delta wc}^s_{ij}\right) + (\overline{B}_{ij} - B_{ij}) \widehat{\delta ws}^s_{ij} \\
&+ G_{ij} \left(\widehat{xv}^s_{ij} - \widehat{wc}^s_{ij}\right) - B_{ij} \widehat{ws}^s_{ij} \quad  \forall ij\in \Edges    
\end{align*}}
\vspace{-0.1cm}
where $\widehat{xv}^s_{ij},\widehat{wc}^s_{ij},\widehat{ws}^s_{ij},\widehat{\delta xv}^s_{ij},\widehat{\delta wc}^s_{ij}$ and $\widehat{\delta ws}^s_{ij}$ satisfy the following constraints:

\vspace{-0.1cm}
{\fontsize{9}{8}\selectfont
\begin{subequations}
\footnotesize
   \label{eq:relaxation}
\begin{align}
   &\widehat{v}^s_{i} \in \langle v^s_i\rangle^{MC-q}, \  \widehat{xv}^s_{ij} \in \langle \tilde{x}^s_{ij}, \widehat{v}^s_{i} \rangle^{MC} \\  
   &\widehat{w}^s_{ij} \in \langle v^s_i, v^s_j \rangle^{MC} \\
   &\widehat{wc}^s_{ij} \in \langle \widehat{w}^s_{ij}, \widehat{cs}^s_{ij} \rangle^{MC}, \
   \widehat{ws}^s_{ij} \in \langle \widehat{w}^s_{ij}, \widehat{sn}^s_{ij} \rangle^{MC} \\ 
   &\widehat{\delta wc}^s_{ij} \in \langle \delta^s_{ij}, \widehat{wc}^s_{ij} \rangle^{MC}, \ 
   \widehat{\delta ws}^s_{ij} \in \langle \delta^s_{ij}, \widehat{ws}^s_{ij} \rangle^{MC} \\ 
 &\widehat{cs}^s_{ij} \leq \tilde{x}^s_{ij} - \left(\frac{1- \cos(\theta^u)}{(\theta^u)^2}\right) ({\theta^s_{ij}}^2+(\tilde{x}^s_{ij}-1)(\theta^{M})^2) \label{eqn:cos_qc}\\
&\tilde{x}^s_{ij}\cos(\theta^u) \leq \widehat{cs}^s_{ij} \leq \tilde{x}^s_{ij}\\
& \widehat{sn}^s_{ij} \leq \cos(\theta^u/2)\theta^s_{ij}+\tilde{x}^s_{ij}(\sin(\theta^u/2)-\theta^u/2\cos(\theta^u/2)) \\ & \nonumber \hspace{20pt}+ (1-\tilde{x}^s_{ij})(\cos(\theta^u/2)\theta^{M} +1)\\
&\widehat{sn}^s_{ij} \geq \cos(\theta^u/2)\theta^s_{ij}-\tilde{x}^s_{ij}(\sin(\theta^u/2)-\theta^u/2\cos(\theta^u/2)) \\ & \nonumber \hspace{20pt}- (1-\tilde{x}^s_{ij})(\cos(\theta^u/2)\theta^{M} +1)\\
&\tilde{x}^s_{ij}\sin(-\theta^u) \leq \widehat{sn}^s_{ij} \leq \tilde{x}^s_{ij}\sin(\theta^u)
\end{align}
\end{subequations}}
Finally, we relax constraint \eqref{eqn:acopf_loss3} into a second-order conic constraint of the form ${p^s}^{2}_{ij} + {q^s}^{2}_{ij} \leq \loss_{ij} {\hat{v}^s}_{i}$.

\vspace{-0.18cm}
\subsection{DC approximation}

Within the expansion planning literature, the DC approximation is often used for designing power systems. For comparison purposes, 
we develop a version of the ORGDT problem based on the DC model in \eqref{eq:master_dc}. 

{\fontsize{9}{8}\selectfont
\begin{subequations}
\begin{flalign}
\SubProblemDC  &  = \{\mathbf{x}^s,\boldsymbol{\tau}^s,\mathbf{t}^s,\boldsymbol{\delta}^s, \boldsymbol{\gamma}^s, \mathbf{zp}^s, \mathbf{u} :  \nonumber \\
& \textbf{On/Off DC power flow equations:} \nonumber \\
& p^s_{ij} \leq -\tilde{B}_{ij} \left(\theta^s_{ij} + \theta^M(1 - \tilde{x}^s_{ij}) \right) \quad  \forall ij\in \Edges,\\
& p^s_{ij} \geq -\tilde{B}_{ij} \left(\theta^s_{ij} - \theta^M(1 - \tilde{x}^s_{ij}) \right) \quad  \forall ij\in \Edges,\\
& \tilde{G}_{ij} = (\bar{G}_{ij} - G_{ij})\delta^s_{ij} + G_{ij} \quad \quad \forall ij \in \Edges \label{eqn:dc_conduct_disj}, \\
& \tilde{B}_{ij} = (\bar{B}_{ij} - B_{ij})\delta^s_{ij} + B_{ij} \quad \quad \forall ij \in \Edges \label{eqn:dc_suscept_disj}, \\
& p^s_{ij} + p^s_{ji} = 0, \ \theta^s_{ij} =  \theta^s_{i} -  \theta^s_{j} + \phi^s_{ij} \quad  \forall ij\in \Edges, \\
& gp^s_{i} - lp^s_{i} = \sum_{ij\in \Edges} p^s_{ij}  \quad \quad \forall i\in \Nodes,  \label{eqn:dc_balance_constraint1}\\
& \textbf{Operational limits and topology constraints:} \nonumber \\
& |p^s_{ij}| \leq \tilde{x}^s_{ij} \capacity_{ij} \quad \quad \forall ij\in \Edges,  \label{eqn:dc_line_constraint_ij} \\  
& |\theta^s_{ij}| \leq \tilde{x}^s_{ij}\theta^u + (1-\tilde{x}^s_{ij})\theta^M \quad \quad  \forall ij \in \Edges, \label{eqn:dc_theta_bound}\\
& -\phi^u \gamma^s_{ij} \leq \phi^s_{ij} \leq \phi^u \gamma^s_{ij} \quad \forall ij\in \Edges,  \label{eqn:dc_phi_bnd} \\
& x^s_{ij} = t^s_{ij} \quad \quad \forall ij\in \Damages_s  \label{eqn:dc_damage}\\
& \tilde{x}^s_{ij} = x^s_{ij} - \tau^s_{ij}\geq 0 \quad \quad \forall ij\in \Edges  \label{eqn:dc_lineminusswitch}\\
& \textbf{Generation and demand constraints} \ \forall i \in \Nodes: \nonumber \\
& 0 \leq zp^s_{i} \leq zp^u_{i}u_{i}, \ 0 \leq gp^s_{i} \leq  gp^u_{i} + zp^s_{i}, \label{eqn:dc_microgrid_constraint_a}\\
& \sum_{i \in \Load} lp^s_{i} \geq lp_{cr} \sum_{i \in \Load} dp_{i},  \label{eqn:dc_critical_constraint1}\\
& \sum_{i \in \Nodes \setminus \Load} lp^s_{i} \geq lp_{ncr} \sum_{i \in \Nodes \setminus \Load} dp_{i} \label{eqn:dc_loadserve_constraint1}\\
& lp^s_{i} =  yp^s_{i} dp_{i}; \ 0 \leq yp^s_i \leq 1, \ \mathbf{x}^s,\boldsymbol{\tau}^s, \mathbf{t}^s \in \{0,1\}  \nonumber \}
\end{flalign}
\label{eq:master_dc}
\end{subequations}
}

\subsection{AC-feasibility formulation}
\label{subsec:ac_feas}

It is important to note that the convex relaxations introduced in section \ref{subsec:convex} 
and the DC approximation of \eqref{eq:master_dc}
might violate the nonconvex AC-equations described in
Eqs. \eqref{eqn:acopf_real_ij}-\eqref{eqn:acopf_reac_ji}, \eqref{eqn:acopf_loss3}. 
In this section, we develop a model for recovering primal feasible solutions to asses the quality of these approximations. This model fixes all the build decisions from the second stage solutions for every disaster $s \in \Disasters$ and solves for the feasibility of the nonconvex equations based on a primal-dual interior point method \cite{wright1997primal}. Formally, let $(\mathbf{\tilde{x}}^s,\mathbf{t}^s,\boldsymbol{\delta}^s, \boldsymbol{\gamma}^s, \mathbf{zp}^s, \mathbf{zq}^s,\mathbf{u})$ be the build decisions obtained in the first stage.
Under the DC model, at bus $i$,
we assume that $zq^s_i$ is equal to half of the real power capacity, that is $|zq^s_i| = zp^s_i/2$. We then minimize a weighted sum of the violations to the resilience constraints which serves as a proxy to measure the distance to AC feasibility.

\vspace{-0.14cm}
{\fontsize{8.5}{8}\selectfont
\begin{subequations}
\begin{flalign}
p_0 := \min \ \ & M\left(\lambda^p_{cr} + \lambda^q_{cr}\right) + \lambda^p_{ncr} + \lambda^q_{ncr}\nonumber \\
\text{s.t.}\quad & \mathrm{Given} \ (\mathbf{\tilde{x}}^s,\mathbf{t}^s,\boldsymbol{\delta}^s, \boldsymbol{\gamma}^s, \mathbf{zp}^s, \mathbf{zq}^s,\mathbf{u}) \nonumber\\ 
&\mathrm{Constraints} \ \eqref{eqn:acopf_real_ij}-\eqref{eqn:ph_ang}, \eqref{eqn:acopf_loss1}-\eqref{eqn:v_bnd}, \eqref{eqn:mg_a}-\eqref{eqn:shed}, \nonumber \\
& \lambda^p_{cr} \geq \left(lp_{cr} \sum_{i \in \Load} dp_{i}\right) - \sum_{i \in \Load} lp^s_{i} , \label{eqn:p_cr_slack}\\
& \lambda^p_{ncr}   \geq \left(lp_{ncr} \sum_{i \in \Nodes \setminus \Load} dp_{i}\right) - \sum_{i \in \Nodes \setminus \Load} lp^s_{i} \label{eqn:p_ncr_slack}\\
& \lambda^q_{cr}  \geq \left(lq_{cr} \sum_{i \in \Load} dq_{i}\right) - \sum_{i \in \Load} lq^s_{i} \label{eqn:q_cr_slack}\\
& \lambda^q_{ncr}  \geq \left(lq_{ncr} \sum_{i \in \Nodes \setminus \Load} dq_{i}\right) - \sum_{i \in \Nodes \setminus \Load} lq^s_{i} \label{eqn:q_ncr_slack}\\
&\lambda^p_{cr}, \lambda^q_{cr}, \lambda^p_{ncr}, \lambda^q_{ncr} \ge 0
\end{flalign}\label{eqn:primal_r}
\end{subequations}}
\vspace{-0.14cm}

\noindent
where $M$ is used to weight satisfaction of critical load and
and $\lambda^p_{cr}, \lambda^q_{cr}, \lambda^p_{ncr}, \lambda^q_{ncr}$ represent the load shed variables on critical and non-critical load constraints (up to the resilience criteria) as shown in constraints \eqref{eqn:p_cr_slack}-\eqref{eqn:q_ncr_slack}.


%% file: Nomenclature.tex

\nomenclature[AA]{$\Nodes$}{set of nodes (buses)}%
\nomenclature[AB]{$\Edges$}{set of edges (lines and transformers) 
}%
\nomenclature[AC]{$\Disasters$}{set of disaster scenarios}%
\nomenclature[AD]{$\Damages_s$}{set of edges that are inoperable during $s\in S$}
\nomenclature[AE]{$\mathbf{i}$}{imaginary number constant}
\nomenclature[AEb]{$\mid.\mid$}{absolute value of the input argument}
\nomenclature[AF]{$c^{x}_{ij}$}{cost to build a line ($i$,$j$) if the line does not already exist}
\nomenclature[AG]{$c^{\tau}_{ij}$}{cost to build a switch on line ($i$,$j$)}
\nomenclature[AH]{$c^t_{ij}$}{cost to harden a line ($i$,$j$)}
\nomenclature[AHb]{$c^{\delta}_{ij}$}{cost to build a FACTS device on line ($i$,$j$)}
\nomenclature[AHc]{$c^{\gamma}_{ij}$}{cost to build a PST device on line ($i$,$j$)}
\nomenclature[AI]{$c^{zp}_{i}$}{cost of real generation capacity at bus $i$}
\nomenclature[AJ]{$c^u_{i}$}{cost to build a generation facility at bus $i$}
\nomenclature[AK]{$T_{ij}$}{apparent power thermal limit on line ($i$,$j$)}
\nomenclature[AL]{$\Load$}{set of buses whose loads are critical}
\nomenclature[AM]{$G_{ij}+\mathbf{i}B_{ij}$}{admittance of line ($i$,$j$)}
\nomenclature[AN]{$R_{ij}+\mathbf{i}X_{ij}$}{impedance of line ($i$,$j$)}
\nomenclature[AM1]{$\bar{G}_{ij}+\mathbf{i}\bar{B}_{ij}$}{modified admittance of line ($i$,$j$) due to top transformer}
\nomenclature[AN1]{$R_{ij}+\mathbf{i}\bar{X}_{ij}$}{modified impedance of line ($i$,$j$) due to top transformer}
\nomenclature[AO]{$\theta^u$}{phase angle difference limit}
\nomenclature[AP]{$\theta^M$}{big-M value on phase angle difference limit}
\nomenclature[AQ]{$\phi^u$}{phase shift limit}
\nomenclature[AR]{$v^l_i,v^u_i$}{lower and upper bound on voltage at bus $i$, respectively}
\nomenclature[AS]{$dp_{i} + \mathbf{i} dq_{i}$}{AC power demand at bus $i$}
\nomenclature[AT]{$gp^u_{i} + \mathbf{i} gq^u_{i}$}{max. existing AC generation capacity at bus $i$}
\nomenclature[AU]{$zp^u_{i} + \mathbf{i} zq^u_{i}$}{max. AC generation capacity that can be built at bus $i$}
\nomenclature[AW]{$lp_{ncr} + \mathbf{i}lq_{ncr}$}{fraction of non-critical AC power loads that must be served}
\nomenclature[AV]{$lp_{cr} + \mathbf{i} lq_{cr}$}{fraction of critical AC power loads that must be served}

\nomenclature[BA]{$x_{ij}$}{determines if line ($i$,$j$) is built}%
\nomenclature[BB]{$\tau_{ij}$}{determines if line ($i$,$j$) has a switch built}
\nomenclature[BC]{$t_{ij}$}{determines if line ($i$,$j$) is hardened}
\nomenclature[BD]{$u_{i}$}{determines the generation capacity built at bus $i$}
\nomenclature[BE]{$\delta_{ij}$}{determines if FACTS device is built on line ($i$,$j$)}
\nomenclature[BF]{$\gamma_{ij}$}{determines if PST device is built on line ($i$,$j$)}
\nomenclature[BG]{$x^s_{ij}$}{determines if line ($i$,$j$) is used during disaster $s \in \cal{S}$}
\nomenclature[BH]{$\tau^s_{ij}$}{determines if switch ($i$,$j$) is used during disaster $s \in \cal{S}$}
\nomenclature[BI]{$t^s_{ij}$}{determines if line ($i$,$j$) is hardened during disaster $s \in \cal{S}$}
\nomenclature[BJ]{$\delta^s_{ij}$}{determines if FACTS device on line ($i$,$j$) is used during disaster $s \in \cal{S}$}
\nomenclature[BK]{$\gamma^s_{ij}$}{determines if PST device on line ($i$,$j$) is used during disaster $s \in \cal{S}$}

\nomenclature[CAa]{$\theta^s_{i}$}{phase angle at bus $i$ during disaster $s \in \cal{S}$}
\nomenclature[CAb]{$\phi^s_{ij}$}{phase angle shift introduced by PST on line $ij$ during disaster $s \in \cal{S}$}
\nomenclature[CB]{$v^s_{i}$}{voltage at bus $i$ during disaster $s \in \cal{S}$}
\nomenclature[CC]{$l^s_{ij}$}{current magnitude squared ($\lvert I_{ij}^s\rvert^2$) on line ($i$,$j$) during disaster $s \in \cal{S}$}
\nomenclature[CD]{$p^s_{ij}+\mathbf{i}q^s_{ij}$}{AC power flow on line ($i$,$j$) during disaster $s \in \cal{S}$}
\nomenclature[CE]{$zp_{i} +\mathbf{i} zq_{i}$}{determines the capacity for AC generation at bus $i$}
\nomenclature[CF]{$zp^s_{i}+\mathbf{i}zq^s_{i}  $}{determines the capacity for AC generation at bus $i$ during disaster $s \in \cal{S}$}
\nomenclature[CG]{$gp^s_{i} + \mathbf{i}gq^s_{i} $}{AC power generated at bus $i$  during disaster $s \in \cal{S}$}
\nomenclature[CH]{$lp^s_{i} + \mathbf{i} lq^s_{i}$}{AC power load delivered at bus $i$ during disaster $s \in \cal{S}$}
\nomenclature[CI]{$yp^s_{i} +\mathbf{i}yq^s_{i}$}{determines the fraction of AC power load served at bus $i$ during disaster $s \in \cal{S}$}


%% file: Algo.tex

\section{Algorithms}
\label{Sec:algo}
In this section we discuss algorithms we use to solve the ORGDT. The ORGDT is a two-stage Mixed-Integer Quadratically Constrained Program (MIQCP) with a block diagonal structure. 
In order to exploit this structure, we generalize the  scenario-based decomposition (SBD) techniques of \cite{yamangil2015resilient} to solve the ORGDT. In the remainder of this paper, let $\ORGDT(\Disasters^\prime)$ denote the ORGDT with the scenario set $\Disasters^\prime \subseteq \Disasters$ and $\sigma$ denote the vector of construction variables $\lineVariable_{ij}$, $\switchVariable_{ij}$, $\hardenVariable_{ij}$,$\gamma_{ij}$, $\delta_{ij}$ for all $ij \in \Edges$ and $\facilityVariable_{i}, zp_i$ for all $i \in \Nodes$. 
SBD can be applied to ORGDT after the following key observation:

\begin{obs}
The second stage variables do not appear in the objective function. Therefore any optimal first stage solution based on a subset of the second stage subproblems that is feasible for the remaining scenarios, is an optimal solution for the original problem.
\end{obs}


The SBD approach is outlined in Algorithm \ref{algo:sbd}. At a high level, Algorithm \ref{algo:sbd} solves problems with iteratively larger sets of scenarios until a solution is obtained that is feasible for all scenarios. The algorithm takes as input the set of scenarios and an initial scenario to consider, $\Disasters^\prime$.  Line 2 solves the ORGDT on $\Disasters^\prime$,
where $\ORGDT(\Disasters^\prime)$ and $\sigma^*$ are used to denote the problem and solution respectively.  Line  3 then evaluates   $\sigma^*$ on the remaining scenarios in $\Disasters \setminus \Disasters^\prime$.
The function $l:P^\prime(s,\sigma^*) \rightarrow \mathbb{R}_+$, is an infeasibility measure that is 0 if the problem is feasible for scenario $s$, positive otherwise. This function is implemented by solving the feasibility subproblem reliability constraints, i.e. total and critical demands are satisfied.
We use a version of (\ref{eqn:primal_r}) where $M=1$ and the non convex constraints are replaced with their convex counterparts.
This function prices the current solution over $s\in \Disasters\setminus \Disasters^\prime$. If all prices are 0, then the algorithm terminates with solution $\sigma^*$ (lines 4-5).  Otherwise, the algorithm adds the scenario with the worst infeasibility measure to $\Disasters^\prime$ (line 7).

\begin{algorithm}[!ht]
\small
\SetKwInOut{Input}{input}
\Input{A set of disasters $\Disasters$ and let $\Disasters^\prime = \Disasters_0$\;}
\While{$\Disasters\setminus \Disasters^\prime \neq \emptyset$}{
	$\sigma^* \leftarrow$ Solve $\ORGDT(\Disasters^\prime)$\;
	$I \leftarrow \left<s_1, s_2 \ldots s_{|\Disasters\setminus \Disasters^\prime|} \right> s\in \Disasters\setminus \Disasters^\prime :$ $l(\ORGDT^\prime(s_i,\sigma^*)) \ge l(\ORGDT^\prime(s_{i+1},\sigma^*))$\;
	\eIf{$l(\ORGDT^\prime(I(0),\sigma^*)) \leq 0$}{
		\Return{$\sigma^*$}\;
	}{
		$\Disasters^\prime \leftarrow \Disasters^\prime \cup  I(0)$\;
	}
\Return{$\sigma^*$}
}
\caption{Scenario-Based Decomposition}
\label{algo:sbd}
\end{algorithm}

\begin{rmk}
We observed that the LP-relaxation for the ORGDT is very loose. To overcome this issue, we augmented every iteration of Algorithm \ref{algo:sbd} with the previous optimum objective value as a lower bound for the current iteration.
\end{rmk}

\noindent
{\bf Cutting-plane algorithm for quadratic constraints:}
Even though optimization theory guarantees that the set of convex inequalities in ORGDT
can be solved efficiently, several numerical experiments demonstrated that it was challenging to solve even moderately sized problems using
state-of-the-art MIQCP solvers (CPLEX/Gurobi). Either the  solver convergence was very slow or it terminated with ``numerical trouble". In order to circumvent this issue, we adopted a cutting-plane approach.
Let $\mathcal{P}^{miq}$ be a general MIQCP which contains the following quadratic constraint (rotated second-order cone):
\begin{equation}
\label{eqn:SOC_gen}
\mathbf{x}^T\mathbf{x} \leq yz, \ y\geq 0,z\geq 0, \ (\mathbf{x},y,z) \in \mathbb{R}^{|\mathbf{x}|+2}
\end{equation}

Algorithm \ref{algo:cp} outlines the cutting plane procedure to solve $\mathcal{P}^{miq}$. The key idea behind this algorithm is to solve a finite sequence of MILPs to obtain an optimal solution for the original MIQCP. Let $\mathcal{P}^{mil}$, a MILP, represent the original problem without the quadratic constraints. The solution in step 3 of Algorithm \ref{algo:cp} is then guaranteed to be a lower bound to $\mathcal{P}^{miq}$. This lower bound is tightened further for every infeasible solution (step 4) by adding a valid cutting plane (step 5). This cutting plane is valid as it is an outer-approximation of the original feasible convex set. This procedure is repeated until the solution obtained is feasible (satisfies step 4), and hence optimal, for $\mathcal{P}^{miq}$.  

\begin{algorithm}[ht]
\small
\caption{Cutting plane algorithm for MIQCPs}
\label{algo:cp}
\textbf{Notation}: $f(\mathbf{x},y) = \frac{\mathbf{x}^T\mathbf{x}}{y}$ \\
\textbf{Input}: {$\mathcal{P}^{mil}$, $tol>0$, $\mathcal{L}=\varnothing$} \\
$\hat{\mathbf{x}},\hat{y},\hat{z}$ $\leftarrow$ Solve $\mathcal{P}^{mil}(\mathbf{x},y,z)$ \\
\While{$\hat{\mathbf{x}}^T\hat{\mathbf{x}} > \hat{y}\hat{z} + tol$}{
	Augment $\mathcal{L}$ with the following cut: 
	$f(\hat{\mathbf{x}},\hat{y}) +\sum_{i=1}^{|\mathbf{x}|} \frac{\partial f(\mathbf{x},y)}{\partial x_i}(x_i-\hat{x}_i) + \frac{\partial f(\mathbf{x},y)}{\partial y}(y-\hat{y}) \leq z$ \\
$\hat{\mathbf{x}},\hat{y},\hat{z}$ $\leftarrow$ Solve $\left(\mathcal{P}^{mil}(\mathbf{x},y,z) \cup \mathcal{L}\right)$
}
\end{algorithm}

%% file: NR.tex
\section{Case Study}
\label{Sec:case}
In this section, we discuss three numerical studies for the ORGDT based on the modified single area IEEE RTS-96 system \cite{wong1999ieee}.
We first discuss the computational benefits associated with using the SBD algorithm. We second compare solutions based on the DC approximation with solutions based on the QC relaxation. Finally, we examine the benefits of FACTS and PST devices. 
All modeling was done using JuMP \cite{dunning2015jump} on a computer with 32 threads, a 2.6GHz Intel 64 bit processor, 25.6MB L3 cache and 64GB of memory. Gurobi 6.5.2 was used for solving the MIQCPs and MILPs (optimality gap 0.1\%) and Knitro 9.1.0 was used for solving the nonlinear programs.

\noindent
\textbf{Test system} We use a modified IEEE single-area RTS-96 system that has 24 buses including 17 load buses, 38 transmission lines and 32 conventional generators \cite{wong1999ieee}. The total installed capacity of the existing generators is 3405 MW. The total load in the system is 2850 MW. We labeled 1740 MW as critical loads. The network was spatially placed in an area of 4250 miles$^2$. Since this test system does not have the co-ordinate data for every bus, we approximately evaluated them based on the lengths of lines \cite{wong1999ieee}.
The admittance, impedance and apparent power thermal limit values on lines are from the standard test case. We assume that FACTS, PSTs and new generation capacities, if chosen, can be built anywhere in the network.  
The remaining parameters are outlined in Table \ref{tab:param}. 

\begin{table}
\caption{Parameters for the test system}
\vspace{-0.2cm}
\footnotesize
\centering
\begin{tabular}{rl}
\toprule
 $c^x_{ij},c^t_{ij},c^{\tau}_{ij}$ & \$1.35m/mile, \$5000/mile, \$1000 \\
 $c^{\delta}_{ij},c^{\gamma}_{ij}$ & \$50000, \$50000\\
 $c^u_i, c^{zp}_i$ & \$0.1m, \$0.817m/MW\\
 $v^l_i,v^u_i,\phi^u,\bar{X}_{ij}$ & 0.95, 1.05, $60\degree$, $0.5X_{ij}$\\
 $lp_{cr},lq_{cr},lp_{ncr},lq_{ncr}$ & 0.99, 0.99, 0.8, 0.8\\
 Total no. of disaster scenarios & 20 \\
 Critical load buses & 1,5,8,11,13,14,15,16,18,19 \\
 New lines which can be built & (1,10),(10,20),(19,24),(6,14),(20,21)\\
 \bottomrule
\end{tabular}
\label{tab:param}
\vspace{-0.4cm}
\end{table}

\noindent
\textbf{Scenario generation}
The damage scenarios are based on line failure probabilities that follow a multivariate Gaussian distribution with a mean placed at the center of the network. A Bernoulli trial was applied to every line for chosen percentiles (\% damage) to generate the random scenarios.
Empirically, we observed that 20 scenarios were sufficient to represent the features of the distribution.

\subsection{Computational Performance of SBD}
\label{subsec:sbd_perf}
Table \ref{tab:times} compares the computational time of SBD 
with the computational time of solving the full model using a MIQCP (Gurobi) solver. Here, we focus on the problem where only new lines, hardened lines, and distributed generations are available as design choices.
On average, SBD is 17x faster when $\theta^u=15\degree$ and 26x faster when $\theta^u=45\degree$.
We highlight in the 80\% damage case that after 20 hours the MIQCP solver could only provide a 7\% optimality gap,
whereas SBD terminated with an optimal solution after 12 minutes. 
It is important to note that when FACTS and PST devices are included as design options, the MIQCP solver fails to find any feasible solution after 10 hours of computation. 

\begin{table}[h]
\caption{Comparison of wall time (sec.) of SBD and MIQCP solvers on problems without FACTS and PST devices. * indicates that the best gap was 7\%. The first column indicates the level of damage. The second/fifth and fourth/seventh columns indicate the wall times of SBD and MIQCP, respectively. The third and sixth columns show the number of iterations for SBD }
\vspace{-0.2cm}
\scriptsize
\centering
\begin{tabular}{crcrrcr}
\toprule
 & \multicolumn{3}{c}{$\theta^u=15\degree$} & \multicolumn{3}{c}{$\theta^u=45\degree$} \\
 \cmidrule(lr){2-4}
 \cmidrule(lr){5-7}
 Damage & SBD & Iters & Full & SBD & Iters & Full \\
 \midrule
 90\% & 634.8 &3 & 5008.2 & 2419.6 & 4 & 13729.2  \\
 80\% & 303.7 & 5 & 1508.7  & 730.1 &5 & 72886.3$^*$ \\
 70\% & 417.7 & 4 & 6566.9 & 485.1 & 4& 2339.0\\ 
 60\% & 188.3 & 2 & 8853.8  & 157.6 & 3 & 3537.5\\
 50\% & 90.2 & 2 & 6156.7  & 53.1 & 1 & 9372.8\\
 \bottomrule
\end{tabular}
\label{tab:times}
\vspace{-0.2cm}
\end{table}

\subsection{AC vs. DC power flow models}
\label{subsec:AcDc}
In this section we scale the loads of the model between 0.8 and 1.25 at .05 intervals. We define the gap between the solutions based on the QC relaxation and solutions based on the DC approximation as

\begin{figure}[htp]
	\centering
	\subfigure[90\% damage, $\theta^u=15^o$]{
	\includegraphics[scale=0.484]{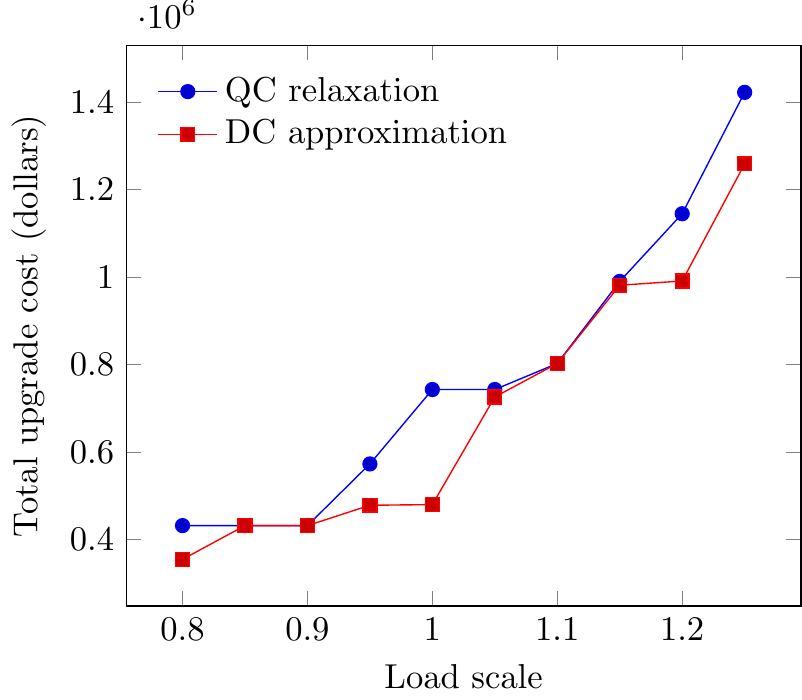}}
	\subfigure[90\% damage, $\theta^u=45^o$]{
	\includegraphics[scale=0.484]{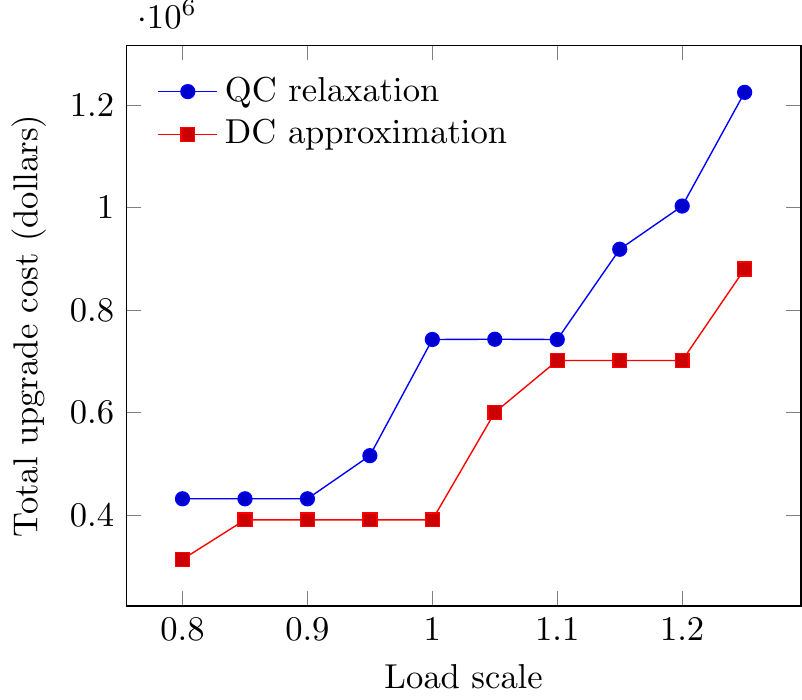}}
	\vspace*{-0.18in}\\
	\caption{Comparison of total upgrade cost for the QC relaxation with the DC approximation}
	\label{fig:gaps}
\end{figure}

{\fontsize{6}{6}\selectfont
$$\zeta = \left(\frac{\mathrm{opt}(\mathcal{P}_0^{AC}) - \mathrm{opt}(\mathcal{P}_0^{DC})}{\mathrm{opt}(\mathcal{P}_0^{AC})}\right)*100$$}
\noindent where opt($\mathcal{P}^{AC}_o$) and opt($\mathcal{P}^{DC}_o$) are optimal upgrade costs for QC relaxation and DC approximation models, respectively.
As shown in table \ref{tab:gaps}, DC-approximation-based solutions are always lower bounds to solutions based on QC relaxations. In the case of small phase angle differences ($\theta^u=15$), the maximum and average $\zeta$ is at the maximum damage case (90\%). When the phase angle differences are increased, the value of $\zeta$ increases dramatically.
This behavior in the solutions is not unexpected. DC approximations tend to perform poorly when there are large phase angle differences as the assumptions behind the approximation no longer hold. 
In Figure \ref{fig:gaps}, we expand the results of the 90\% damage row of Table \ref{tab:gaps} and see how the DC approximation consistently under estimates the design costs.
\vspace{-0.3cm}
\begin{table}[h]
\caption{Each row describes the average, min, and max $\zeta$ across load scaling values between 0.8 and 1.25x for a specified damage level.
}
\vspace{-0.2cm}
\scriptsize 
\centering
\begin{tabular}{ccccccc}
\toprule
 & \multicolumn{3}{c}{$\theta^u=15\degree$} & \multicolumn{3}{c}{$\theta^u=45\degree$} \\
 \cmidrule(lr){2-4}
 \cmidrule(lr){5-7}
 Dam. & avg & min & max & avg & min & max \\
 &gap (\%)&gap (\%)&gap (\%)&gap (\%)&gap (\%)&gap (\%) \\
 \midrule
 90\% & 9.7& 0 & 35.4 & 22.4& 5.5&47.3 \\
 80\% & 9.2& 0 & 25.3& 22.1& 11.7&37.6 \\
 70\% & 7.0& 0& 29.3& 24.0& 15.4&55.7 \\ 
 60\% & 5.6& 0& 29.4& 24.9& 7.4&43.1\\
 50\% & 2.1& 0& 11.3& 35.9& 0&53.1 \\
 \bottomrule
\end{tabular}
\label{tab:gaps}
\end{table}
\vspace{-0.5cm}

\noindent \textbf{AC feasibility analysis} $\zeta$ is not the only metric for comparing the solutions based on the DC and QC models.  It is also important to understand distance to AC feasibility (section \ref{subsec:ac_feas}) in each scenario.
For a given scenario $s$, we calculate the percent apparent-load shed with: 

{\fontsize{3}{3}\selectfont
$$\mu = \left(\frac{\sqrt{(\lambda^p_{cr})^2 + (\lambda^q_{cr})^2}}{\sqrt{\left(\sum_{i\in \mathcal{L}}dp_i\right)^2+\left(\sum_{i\in \mathcal{L}}dq_i\right)^2}}, \frac{\sqrt{(\lambda^p_{ncr})^2 + (\lambda^q_{ncr})^2}}{\sqrt{\left(\sum_{i\in \mathcal{N}\setminus\mathcal{L}}dp_i\right)^2+\left(\sum_{i\in \mathcal{N}\setminus\mathcal{L}}dq_i\right)^2}}\right) \times 100$$}
\noindent for critical and noncritical loads, respectively. 

\begin{figure}[!h]
	\centering
	\subfigure[90\% damage, $\theta^u=15^o$]{
	\includegraphics[scale=0.503]{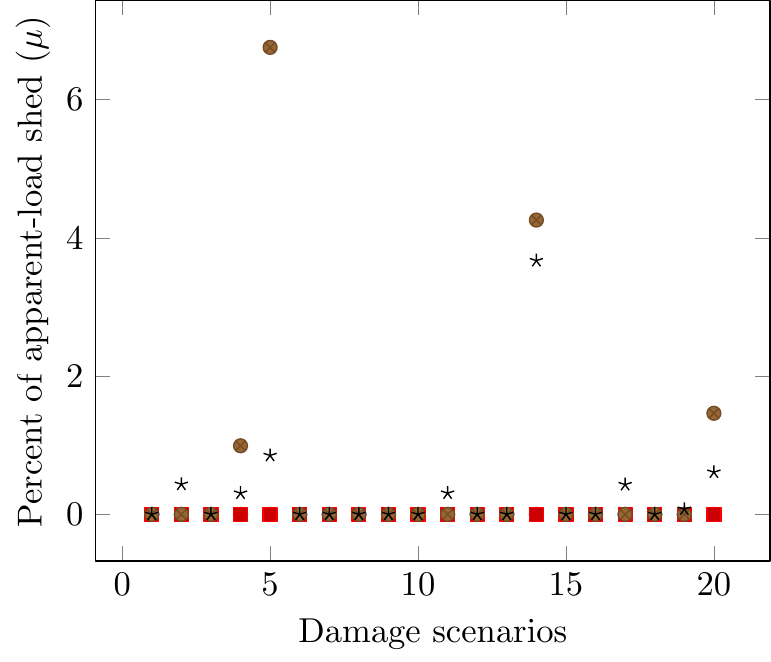}}
	\subfigure[90\% damage, $\theta^u=45^o$]{
	\includegraphics[scale=0.503]{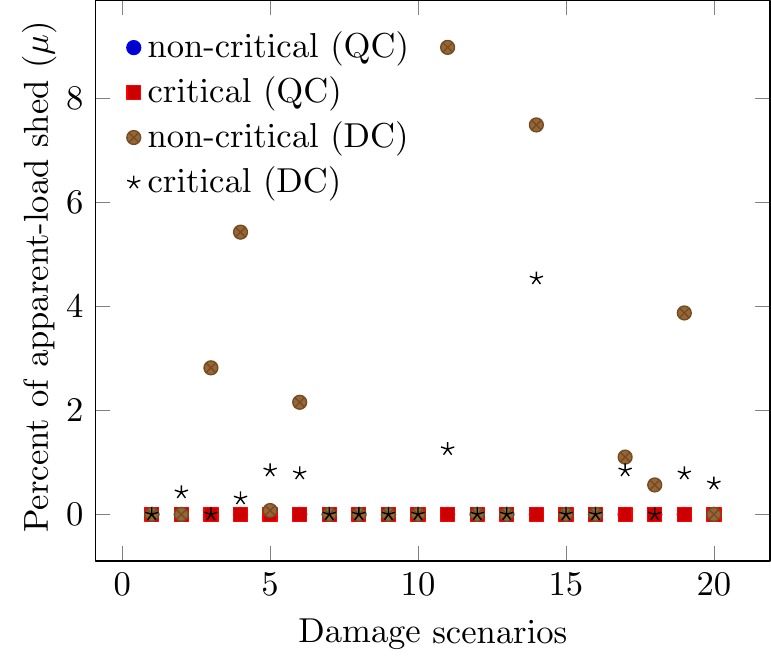}}
	\vspace*{-0.11in}\\
	\caption{Values of $\mu$ for the QC relaxation and DC approximation for 90\% damage scenarios when the load scale is 1.0.}
	\label{fig:shed}
\end{figure}

Table \ref{tab:shed} summarizes the average and maximum $\mu$ for each scenario for the solution found with the DC approximation for a load scale of 1.0.
The QC relaxation is tight and no (additional) load shedding is required. 
Though the DC approximation is considered a good approximation when phase angle differences are small,
we still observe load sheds of up to 3.7\% of critical loads and 6.7\% of non-critical loads in the maximum damage case. Again, the DC approximation's performance degrades further for larger phase angle differences
with a maximum $\mu$ of 4.5\% (critical) and 9.4\% (non-critical). 
Figure \ref{fig:shed} further expands on $\mu$ for each scenario in the 90\% damage case. 

To summarize, OGRDT solutions based on the DC approximation perform reasonable well only on systems with small phase angle differences and limited damage. When the grid is damaged heavily, DC approximation performs poorly, even at smaller $\theta^u$ values. At larger $\theta^u$ values, the DC approximation drastically under approximates the required upgrades, and results in solutions that severely violate resilience criteria based on load shedding.
However, these issues can be addressed by applying AC relaxations in lue of DC approximations. 
\vspace{-0.2cm}
\begin{table}[htp]
\caption{Average and maximum $\mu$ (critical, non-critical) for DC solutions and a load scale of 1.0.}
\vspace{-0.2cm}
\scriptsize 
\centering
\begin{tabular}{cl|l|l|l}
\toprule
 & \multicolumn{2}{c}{$\theta^u=15\degree$} & \multicolumn{2}{c}{$\theta^u=45\degree$} \\
 \cmidrule(lr){2-3}
 \cmidrule(lr){4-5}
 Damage & avg. $\mu$ & max. $\mu$ & avg. $\mu$ & max $\mu$ \\
 \midrule
 90\% & 0.3, 0.7 & 3.7, 6.7 & 0.5, 1.6 &4.5, 8.9  \\
 80\% & 0.1, 1.2& 0.3, 7.3 & 0.2, 3.0& 1.7, 9.2 \\
 70\% & 0.1, 0.1& 0.3, 2.7& 0.15, 2.5& 0.9, 9.4  \\ 
 60\% & 0.1, 0.1& 0.3, 0.2& 0.13, 2.4&0.3, 9.0 \\
 50\% & 0, 0& 0, 0.03& 0, 1.8&0, 7.5 \\
 \bottomrule
\end{tabular}
\label{tab:shed}
\end{table}
\vspace{-0.3cm}

\vspace{-0.5cm}
\subsection{Benefits of FACTS and PST devices}
In this section, we summarize the benefits of including FACTS and PST device options in the ORGDT. To quantify the benefits, we define the following: 

{\fontsize{6}{6}\selectfont
$$\psi = \left(\frac{\mathrm{opt}(\mathcal{P}_0^{AC}) - \mathrm{opt}(\mathcal{P}_0^{AC-dev})}{\mathrm{opt}(\mathcal{P}_0^{AC})}\right)\times 100$$}

\begin{figure}[h]
   \centering
   \includegraphics[scale=0.459]{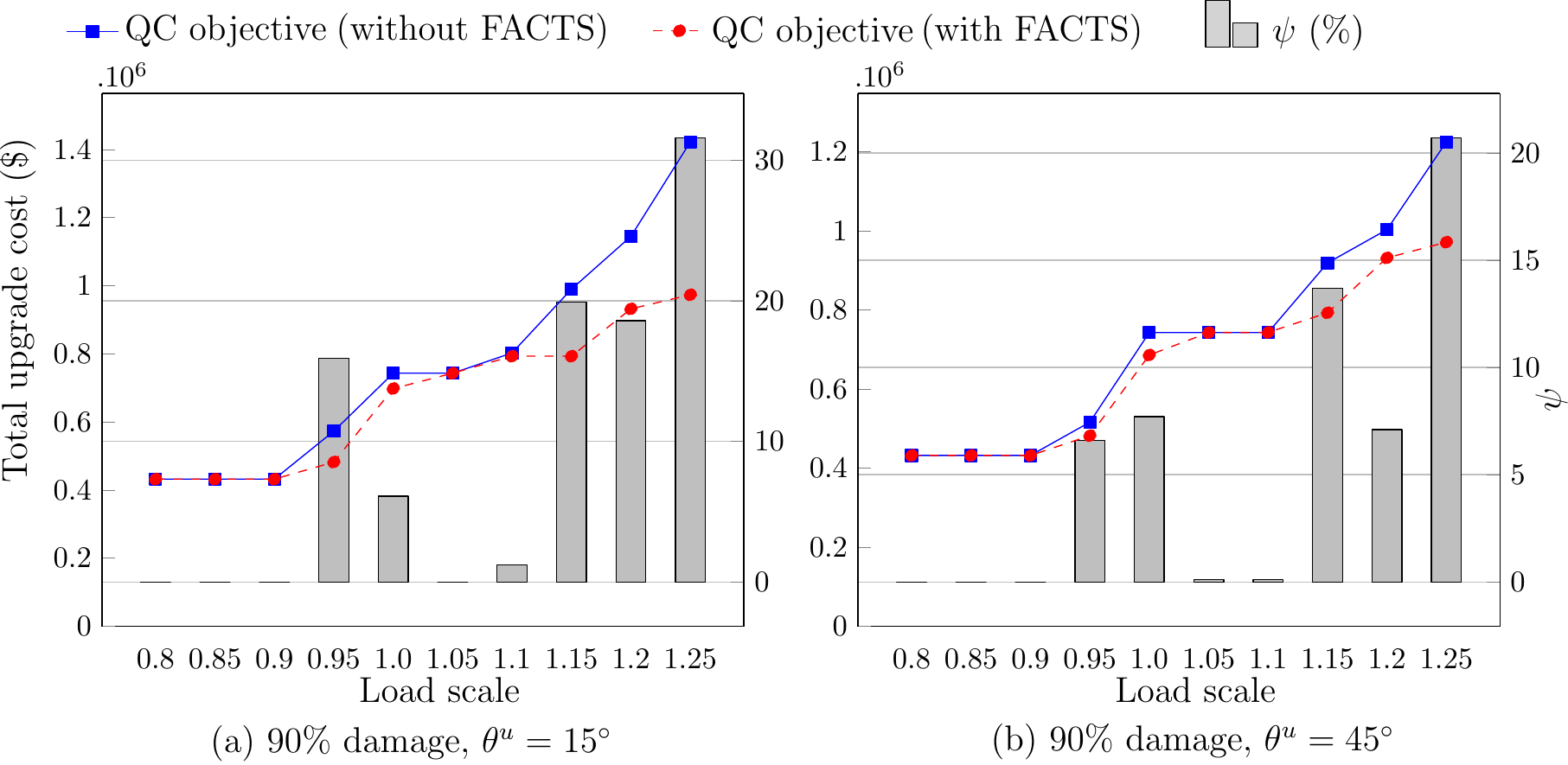}
   	\vspace*{-0.3in}\\
   \caption{Comparison of total upgrade costs and $\psi$ for optimal upgrade solutions with and without FACTS devices.} 
   \label{fig:facts}
\end{figure}

\noindent
\textbf{FACTS devices}
Table \ref{tab:facts} summarizes the benefits of including FACTS devices to achieve advantages of operating transmission grids at larger $\theta^u$ values without the stability issues. Overall, we observe that savings is a much as 17.2\% when $\theta^u=15\degree$.
The benefit of FACTS devices is smaller when $\theta^u$ is larger, since the operating points are less congested and the need for such devices is reduced.
We use Figure \ref{fig:facts} to expand on $\psi$ for the 90\% damage case.

\begin{table}[htp]
\caption{Average and maximum (parenthesis indicate load scale value for the max) $\psi$ across load scales between 0.8 and 1.25 when FACTS devices are included in the OGRDT.
}
\vspace{-0.2cm}
\scriptsize 
\centering
\begin{tabular}{cl|l|l|l}
\toprule
 & \multicolumn{2}{c}{$\theta^u=15\degree$} & \multicolumn{2}{c}{$\theta^u=45\degree$} \\
 \cmidrule(lr){2-3}
 \cmidrule(lr){4-5}
 Damage & avg. $\psi$ & max. $\psi$ & avg.  $\psi$ & max. $\psi$ \\
 \midrule
90\%	&	9.3	&	31.6	(1.25)	&	5.6	&	20.7	(1.25) \\
80\%	&	9.6	&	31.3	(1.2)	&	5.5	&	20.2	(1.25) \\
70\%	&	10.5	&	28.0	(1.2)	&	4.5	&	14.6	(1.25) \\
60\%	&	10.3	&	29.7	(1.2)	&	4.6	&	28.0	(1.25) \\
50\%	&	17.2	&	33.2	(1.2)	&	11.8	&	29.4	(1.25) \\
 \bottomrule
\end{tabular}
\label{tab:facts}
\end{table}

\vspace{-0.2cm}
In Figure \ref{fig:facts}, it is interesting to note the (sometimes) non-increasing trend in $\psi$ as load scale increases.
Intuitively, this indicates that upgrades required for certain load scales are also sufficient for further increases in load. More importantly, perhaps, the OGRDT solution with and without FACTS devices is sometimes identical. This phenomena indicates that at certain load scales, FACTS devices cannot eliminate the need for other upgrades (line hardening, new lines, etc.) and once those upgrades are installed, the need for FACTS devices is eliminated. Of course, as the load scale increases, FACTS devices are once again required.  


\noindent
\textbf{PST devices} Generally speaking, we only found PST devices to be of benefit when the network is congested and has larger load scales.
Under the 90\% damage cases, $\theta^u=15\degree$ and load scales greater than 1.20, we found solutions with savings of around 10\%. Though the improvement in $\psi$ when PSTs are included is not as impressive as when FACTS devices are included, we believe that the savings are non trival when considering large-scale grids.

%% file: Conclusions.tex


\vspace{-.2cm}

\section{Conclusions}
\label{Sec:conc}

In this paper, we formulated, modeled, and developed a model of ORGDT. 
Our contributions include a comparison of solutions to the ORGDT with the DC approximation and the QC relaxation that demonstrates the necessity of AC power flow modeling in resilience studies. We also developed an approach for recovering fully AC feasible solutions that indicate that solutions based on the QC relaxation tend to be tight in practice. Finally, we have also shown that a small set of strategically located PST and FACTS devices as design options can reduce the need for other design considerations like line hardening or distributed generation.


There remain a number of interesting future directions for work in this area. First, it will be important to scale the approaches of the paper to larger, more realistic transmission grids. An important idea in this area is to limit the upgrades to those that are deemed practical by subject matter experts. Second, it will be important to introduce tighter and improved models of FACTS and PST devices (such as continuous set points). Finally, it will be interesting to consider= approaches for using FACTS and PST devices to improve the restoration process of transmission grids, as seen in \cite{Mak2014}.


